\documentclass[preprint,12pt]{elsarticle}

\usepackage{amsmath,amssymb,amsthm,mathtools}
\usepackage{microtype}
\usepackage[T1]{fontenc}
\usepackage{lmodern}
\usepackage{xcolor}

\newcommand{\rev}[1]{\textcolor{red}{#1}}
\newenvironment{revblock}{\begingroup\color{red}}{\endgroup}

\journal{Fuzzy Sets and Systems}

\newtheorem{theorem}{Theorem}[section]
\newtheorem{lemma}[theorem]{Lemma}
\newtheorem{proposition}[theorem]{Proposition}
\newtheorem{corollary}[theorem]{Corollary}
\theoremstyle{definition}
\newtheorem{definition}[theorem]{Definition}
\newtheorem{example}[theorem]{Example}
\theoremstyle{remark}
\newtheorem{remark}[theorem]{Remark}

\newcommand{\Der}{\operatorname{Der}}
\newcommand{\PDer}{\operatorname{PDer}}
\newcommand{\Fix}{\operatorname{Fix}}
\newcommand{\id}{\operatorname{id}}
\newcommand{\pt}{\mathrm{pt}}
\newcommand{\TP}{T_{\mathrm P}}
\newcommand{\TL}{T_{\mathrm L}}
\newcommand{\TG}{T_{\mathrm G}}

\begin{document}

\begin{frontmatter}

\title{\rev{Leibniz Derivations on Ordered Incline-Type Algebras: Distributivity Detection and Continuous t-Norm Classification}}

\author{Qiuyu Shi}
\author{Wei Yao}
\author{Yi Zhang}
\address{School of Mathematics and Statistics, Nanjing University of Information Science and Technology, Nanjing 210044, China}
\ead{corresponding@email.example}

\begin{abstract}
\begin{revblock}
We study Leibniz-type derivations on a weak class of commutative ordered algebras in which the additive structure is a bounded join-semilattice and multiplication is an isotone commutative monoid operation, while distributivity over joins is not assumed. Every multiplication map $p_a(x)=ax$ is a derivation, and the assignment $a\mapsto p_a$ together with evaluation at the unit forms a canonical Galois insertion. This principal representation allows distributivity and residuation to be detected rather than assumed: pointwise binary-join closure of the derivation space characterizes inclines, and in the complete case arbitrary pointwise-join closure characterizes commutative integral quantales. Principal spectra further encode idempotence, meet behavior, finite-join preservation, and residuation.

For commutative integral quantales, the full poset of Leibniz derivations is a complete lattice, although individual derivations are not required to preserve arbitrary joins. Evaluation at the unit is a surjective complete-lattice homomorphism whose fibers are complete intervals with canonical principal minima. For every t-norm, endpoint-modified principal maps form a natural family of derivations; for continuous Archimedean t-norms they exhaust all derivations. Consequently, all continuous Archimedean t-norms have the same fibered derivation poset, so order and evaluation fibers forget the strict/nilpotent type, whereas composition of principal derivations reconstructs the multiplication. For general continuous t-norms, canonical cut maps detect exactly the idempotent skeleton, and finite composition powers of principal derivations identify the Product- and {\L}ukasiewicz-type ordinal-sum components.
\end{revblock}
\end{abstract}


\end{frontmatter}

\section{Introduction and preliminaries}

\begin{revblock}
Derivations have been transported from ring theory to a broad range of ordered and many-valued algebraic structures. The lattice-theoretic study goes back to Sz\'asz \cite{Szasz1975} and was developed further by Ferrari \cite{Ferrari2001} and by Xin, Li and Lu, Xin, and Gan and Guo \cite{XinLiLu2008,Xin2012,GanGuo2022}. For the classical lattice notion, Ferrari proved in particular that pointwise joins and meets of derivations are again derivations in the distributive case. Our definition is weaker: only the Leibniz identity is imposed, with no a priori join-preservation requirement on the derivation itself.

In incline theory, Alshehri considered derivations satisfying $d(xy)=d(x)y+x d(y)$, where addition is idempotent \cite{AlshehriIncline2010}. Subsequent work developed right and generalized derivations of incline algebras \cite{KimRight2013,OzbalFirat2014}. Those studies take place inside incline-type structures in which distributivity is already part of the ambient axioms. One aim of the present paper is instead to start before distributivity is known and to recover it from the pointwise-join behavior of the derivation space.

Derivations have also been investigated in algebras related to fuzzy logic. Torkzadeh and Abbasian introduced $(\odot,\vee)$-derivations for BL-algebras \cite{TorkzadehAbbasian2013}, and Xin, Feng and Yang obtained derivation-based characterizations of G\"odel algebras \cite{XinFengYang2016}. For MV-algebras, see Alshehri \cite{AlshehriMV2010} and Zhao, Gan and Yang \cite{ZhaoGanYang2024}; structural representations through other types of derivations and Galois connections were studied by Wang, She and Qian \cite{WangSheQian2019}. The same Leibniz pattern has recently been extended to EMV- and EBL-algebras \cite{LiuLiu2024}. In particular, Liu and Liu also obtained a threshold cut-map criterion for linearly ordered EBL-algebras; Section~\ref{sec:cut-spectrum} below identifies the exact necessary and sufficient threshold condition in the continuous t-norm setting. Higher-order variants were developed by Zhao and Yang \cite{ZhaoYang2025}. Related theories occur in GMV-algebras and generalized algebras of BCK logic \cite{RachunekSalounova2018,Ciungu2021}.

The work of Zhao, Gan and Yang is particularly relevant to our later {\L}ukasiewicz analysis. They construct endpoint modifications of principal derivations and prove that such maps exhaust all derivations on finite MV-chains \cite{ZhaoGanYang2024}. Our t-norm specialization addresses the continuum: we prove that the same form exhausts all derivations of the standard {\L}ukasiewicz chain $[0,1]$, obtain an independent Product classification, and transport the two model classifications to every continuous Archimedean t-norm.

The residuated-lattice work of He, Xin and Zhan is also close to the present setting. They studied multiplicative derivations, principal ideal derivations, their adjoints and fixed-point sets, and applications to Heyting algebras \cite{HeXinZhan2016}; subsequent developments include \cite{Kondo2018,Kondo2020}. However, in a residuated lattice each multiplication map is already a left adjoint and therefore preserves all existing joins. Likewise, in an incline, distributivity is part of the axioms. Consequently, neither distributivity nor residuation can be detected once one starts inside those categories.

Principal derivations have also been used as reconstruction data in multiplicative semilattice theory. Wang, Borumand Saeid and Wang showed that an idempotent commutative quantale is determined by its principal derivations \cite{WangBorumandWang2018}. Their result starts with a distributive multiplicative structure. By contrast, Theorems~\ref{thm:binary-join} and \ref{thm:quantale-characterization} below begin before distributivity is assumed and characterize exactly when pointwise joins of derivations force the ambient algebra to become an incline or a quantale.

A related distinction is important in the quantale literature. The join-preserving notion of quantale derivation was introduced by Xiao and Li \cite{XiaoLi2012} and developed further by Xiao and Liu \cite{XiaoLiu2016}. Here a derivation is required only to satisfy the Leibniz identity. Thus our complete-lattice result does not assume arbitrary-join preservation of individual derivations; pointwise joins become derivations because the underlying multiplication is quantalic.

The organizing observation is that, for every element $a$, the multiplication operator $p_a(x)=ax$ is a derivation even before distributivity is assumed. Moreover, with $\varepsilon(d)=d(1)$, one has $p_a\le d$ if and only if $a\le\varepsilon(d)$. Hence $a\mapsto p_a$ and evaluation at the unit form a Galois insertion, and $p_{d(1)}$ is the largest principal derivation below $d$.

The paper develops three groups of results. First, pointwise behavior of derivations detects algebraic laws. We prove that $K$ is an incline if and only if $\Der(K)$ is closed under pointwise binary joins; in the complete case, arbitrary pointwise-join closure characterizes commutative integral quantales. Principal spectra then separate idempotence, meet behavior, finite-join preservation, and residuation.

Second, for a commutative integral quantale, $\Der(K)$ is a complete lattice. Joins are pointwise, while meets are obtained by applying a canonical derivation interior to pointwise meets. Evaluation at the unit fits into an adjoint triple $\eta\dashv\varepsilon\dashv u$, and every fiber $\Der_a(K)$ is the complete interval $[p_a,u_a]$.

Third, we specialize to t-norms. Composition of principal derivations reconstructs the ordered t-norm monoid. We show that endpoint-modified principal maps are derivations for every t-norm and that, for a continuous Archimedean t-norm, they are all the derivations. Consequently all continuous Archimedean t-norms share the same fibered derivation poset, so this invariant does not retain the strict/nilpotent type. For general continuous t-norms, canonical cut maps detect the idempotent skeleton, while finite principal composition powers identify the types of the Archimedean ordinal-sum components.

The paper is organized as follows. Section~2 introduces ordered incline-type algebras and Leibniz derivations. Section~3 develops the principal representation and characterizes distributivity by pointwise joins. Section~4 studies the principal spectra. Section~5 establishes the complete derivation lattice and its fiber decomposition. Section~6 treats t-norm reconstruction and the continuous Archimedean classification. Section~7 identifies the cut spectrum with the idempotent skeleton. Section~8 gives an operator form of the strict/nilpotent power criterion and illustrates it on a mixed ordinal sum. Section~9 concludes.
\end{revblock}

\section{Ordered incline-type algebras}

\begin{definition}
An \emph{ordered incline-type algebra} is an algebra
$K=(K,\vee,\cdot,0,1)$
such that:
\begin{enumerate}
\item $(K,\vee,0,1)$ is a bounded join-semilattice;
\item $(K,\cdot,1)$ is a commutative monoid;
\item the multiplication $\cdot$ is isotone in each variable w.r.t. the natural order.
\end{enumerate}
Notice that, there is no distributivity of multiplication over $\vee$ assumed.
\end{definition}

We usually omit the multiplication dot and write $xy$ instead of $x\cdot y$. It is easily seen that for an ordered incline-type algebra $K$ is and that $x,y\in K$, we have $xy\le x$, $xy\le y$ and $x0=0$.

\begin{remark}
If $x(y\vee z)=xy\vee xz\ (\forall x,y,z\in K)$, then $K$ is a bounded commutative incline in the usual sense. Inclines and incline matrices are surveyed in \cite{KimRoush2004,CaoKimRoush1984}.
\end{remark}

\begin{example}
The following are ordered incline-type algebras:
\begin{enumerate}
\item every bounded commutative incline algebra is an ordered incline-type algebra.
	
\item every bounded lattice $L$ with multiplication $x\cdot y=x\wedge y$ is an ordered incline-type algebra. 

\item every residuated lattice, hence every integral commutative quantale,  is an ordered incline-type algebra. 

\item the unit interval $[0,1]$ with a t-norm is an ordered incline-type algebra.
\end{enumerate}
\end{example}

\begin{definition}
For an ordered incline-type algebra $K$, a mapping $d:K\longrightarrow K$ is called a \emph{derivation} on $K$ if
$$ d(xy)=d(x)y\vee x d(y)\ (\forall x,y\in K).$$
We write $\Der(K)$ for the set of all derivations on $K$. 
\end{definition}

This definition specializes to the usual lattice Leibniz rule when multiplication is meet, to the incline derivation rule of \cite{AlshehriIncline2010}, and to the multiplicative derivation considered in residuated lattices in \cite{HeXinZhan2016}. It is weaker than the quantale notion used in \cite{XiaoLiu2016}, where arbitrary-join preservation is imposed as an additional axiom.

\begin{lemma}\label{lem:d-basic}
Suppose that $K$ is an ordered incline-type algebra and that $d:K\longrightarrow K$ is a derivation. Then the following conditions hold:
\begin{enumerate}
\item $d(0)=0$;
\item $x d(1)\le d(x)\ (\forall x\in K)$;
\item $d(x^2)=x d(x)\ (\forall x\in K)$.
\end{enumerate}
\end{lemma}

\begin{proof}
(1) Let $x=y=0$, we have $d(0)=d(0)0\vee0d(0)=0$.

(2) Since $x=x1$, we have $d(x)=d(x)1\vee x d(1)=d(x)\vee x d(1)$ and then $x d(1)\le d(x)$.

(3) Putting $y=x$, we have $d(x^2)=d(x)x\vee x d(x)=x d(x)$.
\end{proof}

\begin{definition}
Let $(P,\leq)$ be a partially ordered set. A mapping
$f:P\to P$ is called an \emph{interior operator} if the following conditions hold for all $x,y\in P$:
\begin{enumerate}
\item $x\leq y \implies f(x)\leq f(y);$
\item $f(x) \leq x;$
\item $f(f(x))=f(x).$
\end{enumerate}
\end{definition}

\begin{definition}
Let $A$ and $B$ be structures equipped with binary operations $\vee$ and $\cdot$. A mapping $f:A\to B$ is called an \emph{isomorphism} if it is bijective and, for all $a,b\in A$,
\[
f(a\vee b)=f(a)\vee f(b),
\qquad
f(a\cdot b)=f(a)\cdot f(b).
\]
\end{definition} 

\begin{definition}
Let $(M,\leq_M,\cdot_M,1_M)$ and $(N,\leq_N,\cdot_N,1_N)$ be partially ordered monoids. A mapping $f:M\to N$ is called an \emph{order-monoid isomorphism} if it is bijective and, for all $x,y\in M$,
\[
x\leq_M y \iff f(x)\leq_N f(y),\qquad f(x\cdot_M y)=f(x)\cdot_N f(y),\qquad f(1_M)=1_N.
\]
\end{definition}

\section{Principal derivations and distributivity}

For $a\in K$, define $$p_a:K\longrightarrow K,\ p_a(x)=ax\ (\forall x\in K).$$ It is easy to verify that every $p_a$ is a derivation, and so we called it the principal derivation of $a$. We write $\PDer(K):=\{p_a:a\in K\}$ for the set of all principal derivations.

\begin{theorem}\label{thm:principal}
Suppose that $K$ is an ordered incline-type algebra. Then
\begin{enumerate}

\item Define $\varepsilon(d)=d(1)\ (\forall d\in{\rm Der}(K))$. Then $(\eta,\varepsilon)$ forms a Galois adjunction from $K$ to ${\rm Der}(K)$, that is, for all $a\in K$ and all $d\in{\rm Der}(K)$, it holds that 
\[
\eta(a)\le d\Longleftrightarrow a\le \varepsilon(d).
\]

\item $\varepsilon\eta=\id_K$, $\varepsilon$ is a surjection, and $\eta$ is an order-monoid isomorphism; in particular $K\cong {\rm PDer}(K)$.

\item define $d\sim e$ by $d(1)=e(1)$. Then there are an order-monoid isomorphism
\[
\Der(K)/{\sim}\ \cong {\rm PDer}(K).
\]

\item $P(d):=p_{d(1)}$ is the largest principal derivation below $d$, $P$ is an interior operator on $\Der(K)$, and $\Fix(P)=\PDer(K)$.
\end{enumerate}
\end{theorem}

\begin{proof}
(1) If $\eta(a)=p_a\le d$, then $a=p_a(1)\le d(1)=\varepsilon(d)$. Conversely, if $a\le \varepsilon(d) = d(1)$, then by isotonicity and Lemma~\ref{lem:d-basic}(2), $ax\le d(1)x\le d(x)$ for all $x\in K$, so $p_a\le d$.

(2) $\varepsilon\eta(a)=p_a(1)=a$, so $\varepsilon\eta=\id_K$ and $\varepsilon$ is surjective. Then $\eta$ is injective, hence a bijection from $K$ onto $\PDer(K)$. If $a\le b$, then $ax\le bx$ for all $x$, so $p_a\le p_b$; conversely, $p_a\le p_b$ implies $a=p_a(1)\le p_b(1)=b$. By associativity and commutativity, $(p_a\circ p_b)(x)=a(bx)=(ab)x=p_{ab}(x)$, and $\eta(1)=p_1=\id$. Hence $\eta$ is an order-monoid isomorphism and $K\cong\PDer(K)$.

(3) The two sets have the same image under $\varepsilon$, namely $\varepsilon(\Der(K))=K=\varepsilon(\PDer(K))$. Since $\varepsilon(p_a)=a$, the restriction $\varepsilon|_{\PDer(K)}$ is an order-monoid isomorphism onto $K$, with inverse $\eta$ in (2). The kernel equivalence relation of $\varepsilon$ on $\Der(K)$ is exactly $\sim$, so the canonical factorization of the surjection $\varepsilon$ gives $\Der(K)/{\sim}\cong K$, where the quotient carries the order and multiplication transported from $K$: $[d]\le[e]\iff d(1)\le e(1)$ and $[d][e]:=[p_{d(1)e(1)}]$, with unit $[p_1]$. Identifying both sides with the common image $K$, the map $[d]\mapsto p_{d(1)}$ is an order-monoid isomorphism $\Der(K)/{\sim}\cong\PDer(K)$.

(4) By (1), $p_c\le d$ if and only if $c\le d(1)$, so $P(d)=p_{d(1)}$ is the largest principal derivation below $d$. If $d\le e$, then $d(1)\le e(1)$ and $P(d)\le P(e)$; moreover $P(d)\le d$, and $P(d)(1)=d(1)$ implies $P(P(d))=P(d)$. Hence $P$ is an interior operator. Finally, $P(p_a)=p_a$, and $P(d)=d$ forces $d=p_{d(1)}\in\PDer(K)$, so $\Fix(P)=\PDer(K)$.
\end{proof}

\begin{theorem}
${\rm PDer}(K)$ is an ordered commutative incline-type algebra.
\end{theorem}

\begin{proof}
$\PDer(K)$ is ordered pointwise, with composition as multiplication. By Theorem~\ref{thm:principal}(2), $\eta:a\mapsto p_a$ is an order-monoid isomorphism from $K$ onto $\PDer(K)$. An order isomorphism preserves existing joins: since $a\vee b$ exists in $K$, $p_{a\vee b}$ is the least upper bound of $p_a$ and $p_b$ in $\PDer(K)$, and $p_0$, $p_1=\id$ are the least and greatest elements. Hence $\PDer(K)$ is a bounded join-semilattice with $p_a\vee p_b=p_{a\vee b}$. Since $\eta$ preserves multiplication and the unit, $(\PDer(K),\circ,p_1)$ is a commutative monoid.By isotonicity of composition , $K$ have the order isomorphism. Thus $\PDer(K)$ is an ordered incline-type algebra.
\end{proof}

The first structural consequence is that distributivity can be recovered from pointwise joins of derivations.

\begin{revblock}
Ferrari proved pointwise-join closure for derivations of a distributive lattice \cite{Ferrari2001}. The next theorem supplies the converse in the present weaker Leibniz-only setting: closure of the derivation space itself forces distributivity.
\end{revblock}

\begin{theorem}\label{thm:binary-join}
Suppose that $K$ is an ordered incline-type algebra. Then the following conditions are equivalent:
\begin{enumerate}
\item $K$ is an incline, i.e.
 $x(y\vee z)=xy\vee xz\ (\forall x,y,z\in K)$;
\item $\PDer(K)$ is closed under pointwise binary joins;
\item $ p_{a\vee b}=p_a\vee_{\pt}p_b\ (\forall a,b\in K)$;
\item $\Der(K)$ is closed under pointwise binary joins.
\end{enumerate}
\end{theorem}

\begin{proof}
(1)$\Longrightarrow$(4) let $d,e\in\Der(K)$ and define $h:=d\vee_{\pt}e$. By distributivity,
\begin{align*}
 h(xy)
 &=d(xy)\vee e(xy)\\
 &=d(x)y\vee xd(y)\vee e(x)y\vee xe(y)\\
 &=(d(x)\vee e(x))y\vee x(d(y)\vee e(y))\\
 &=h(x)y\vee xh(y).
\end{align*}
Thus $h$ is a derivation.

(4)$\Longrightarrow$(3) fix $a,b\in K$. Then $h:=p_a\vee_{\pt}p_b$ is a derivation and $h(1)=a\vee b$. By Lemma~\ref{lem:d-basic}, $x(a\vee b)=xh(1)\le h(x)=ax\vee bx$.
Since multiplication is isotone and $a,b\le a\vee b$,we have $ax \vee bx \leq x(a \vee b)$. Hence $h(x)=x(a\vee b)=p_{a\vee b}(x)$ for every $x$.

(3)$\Longrightarrow$(2) if $p_{a\vee b}=p_a\vee_{\pt}p_b$, then $p_a\vee_{\pt}p_b$ is principal for all $a,b\in K$.

(2)$\Longrightarrow$(1) if $\PDer(K)$ is closed under pointwise binary joins, then $h=p_a\vee_{\pt}p_b$ is principal, say $h=p_c$. Let $x=1$ then $p_c(1)=p_a(1)\vee_{\pt}p_b(1)$, $i.e.\ c=a\vee b$, so $h=p_{a\vee b}$. Therefore $(a\vee b)x=p_{a\vee b}(x)=p_a(x)\vee p_b(x)=ax\vee bx$, and by commutativity we can get the incline distributivity law.
\end{proof}

\begin{remark}
Even without distributivity, the principal representation identifies $p_{a\vee b}$ as the least principal derivation above both $p_a$ and $p_b$. Theorem~\ref{thm:binary-join} says that incline distributivity is exactly the condition under which this intrinsic principal join agrees with the pointwise join.
\end{remark}

The complete version detects quantale distributivity.

\begin{theorem}\label{thm:quantale-characterization}
Suppose that $K$ is an ordered incline-type algebra whose underlying ordered set is a complete lattice. Then the following conditions are equivalent:
\begin{enumerate}
\item $K$ is a commutative integral quantale, i.e.
\[
 x\Bigl(\bigvee_{i\in I}a_i\Bigr)=\bigvee_{i\in I}xa_i,
\]
for every family $\{a_i\}_{i\in I}\subseteq K$;
\item $\PDer(K)$ is closed under arbitrary pointwise joins;
\item $\Der(K)$ is closed under arbitrary pointwise joins.
\end{enumerate}
When these conditions hold,
\[
 \left(\bigvee_{i\in I}^{\Der}d_i\right)(x)
 =\bigvee_{i\in I}d_i(x),
\]
for every family of derivations.
\end{theorem}

\begin{proof}
(1)$\Longrightarrow$(3) and (1)$\Longrightarrow$(2) let $\{d_i\}_{i\in I}\subseteq\Der(K)$. Put $d=\bigvee_{i\in I}^{\pt}d_i$.
By arbitrary distributivity of multiplication,
\begin{align*}
 d(xy)
 &=\bigvee_i d_i(xy)\\
 &=\bigvee_i\bigl(d_i(x)y\vee xd_i(y)\bigr)\\
 &=\left(\bigvee_i d_i(x)\right)y
   \vee x\left(\bigvee_i d_i(y)\right)\\
 &=d(x)y\vee xd(y).
\end{align*}
Moreover, $\left(\bigvee_i^{\pt}p_{a_i}\right)(x)=\bigvee_i a_i x=\left(\bigvee_i a_i\right)x=p_{\bigvee_i a_i}(x)$.

(3)$\Longrightarrow$(1) let $\{a_i\}_{i\in I}\subseteq K$. The pointwise join $h:=\bigvee_i^{\pt}p_{a_i}$ is a derivation and $h(1)=\bigvee_i a_i$.
By Lemma~\ref{lem:d-basic}, $x\left(\bigvee_i a_i\right)=xh(1)\le h(x)=\bigvee_i xa_i$.
By isotonicity we get the opposite inequality. Hence multiplication preserves arbitrary joins.

(2)$\Longrightarrow$(1) for $h=\bigvee_i^{\pt}p_{a_i}$ there is some $c\in K$ with $h=p_c$. When it takes the value 1 then we have $c=\bigvee_i a_i$. Therefore $x\left(\bigvee_i a_i\right)=p_c(x)=h(x)=\bigvee_i xa_i$.
When (1) holds, the displayed pointwise formula for joins in $\Der(K)$ follows from the first part of the proof.
\end{proof}

\section{Principal-derivation spectra}

The principal representation separates several algebraic phenomena that are automatically conflated in more restrictive ambient categories.

\begin{definition}
For an ordered incline-type algebra $K$, define
\begin{align*}
 E(K)&=\{a\in K:p_a^2=p_a\},\\
 M(K)&=\{a\in K:\Fix(p_a)=\downarrow a\},\\
 D(K)&=\{a\in K:p_a\text{ preserves binary joins}\},\\
 R(K)&=\{a\in K:p_a\text{ has a right adjoint}\}.
\end{align*}
We call these the \emph{idempotent}, \emph{meet}, \emph{distributive}, and \emph{residuated principal spectra}, respectively.
\end{definition}

\begin{proposition}\label{prop:spectra-basic}
Suppose that $K$ is an ordered incline-type algebra. Then:
\begin{enumerate}
\item $E(K)=\{a\in K:a^2=a\}$;
\item $M(K)\subseteq E(K)$ and $R(K)\subseteq D(K)$;
\item $E(K)$ is closed under both $\vee$ and multiplication;
\item $M(K)$, $D(K)$, and $R(K)$ are closed under multiplication.
\end{enumerate}
\end{proposition}

\begin{proof}
(1) For every $a,x\in K$,
\[
 p_a^2(x)=p_a(ax)=a(ax)=a^2x=p_{a^2}(x),
\]
so $p_a^2=p_{a^2}$. Since $a\mapsto p_a$ is injective by Theorem~\ref{thm:principal}(3), $p_a^2=p_a$ is equivalent to $a^2=a$. Hence $E(K)=\{a\in K:a^2=a\}$.

(2) If $a\in M(K)$, then $a\in\downarrow a=\Fix(p_a)$, so $p_a(a)=a$, i.e. $a^2=a$. Thus $M(K)\subseteq E(K)$. If $a\in R(K)$, then $p_a$ has a right adjoint and is therefore a left adjoint, hence preserves all existing joins and in particular binary joins. Thus $R(K)\subseteq D(K)$.

(3) For the closure of $E(K)$ under multiplication, let $e,f\in E(K)$. Then
\[
 (ef)^2=e^2f^2=ef,
\]
so $ef\in E(K)$. For the closure under $\vee$, again let $e,f\in E(K)$. Since $e=e^2\le(e\vee f)^2$ and $f=f^2\le(e\vee f)^2$, we have $e\vee f\le(e\vee f)^2$. On the other hand, by Lemma~\ref{lem:basic-product}, $(e\vee f)^2\le e\vee f$. Hence $(e\vee f)^2=e\vee f$, i.e. $e\vee f\in E(K)$.

(4) For $a,b\in M(K)$, we show $\Fix(p_{ab})=\downarrow ab$. If $x\le ab$, then by Lemma~\ref{lem:basic-product} we have $x\le ab\le a,b$, and since $a,b\in M(K)$, $ax=x=bx$. Hence
\[
 (ab)x=a(bx)=ax=x,
\]
so $x\in\Fix(p_{ab})$. Conversely, if $(ab)x=x$, then by Lemma~\ref{lem:basic-product}, $x=(ab)x\le ab$, so $x\in\downarrow ab$. Thus $ab\in M(K)$.

For $a,b\in D(K)$, we use $p_{ab}=p_a\circ p_b$. Since both $p_a$ and $p_b$ preserve binary joins,
\begin{align*}
 p_{ab}(x\vee y)
 &=p_a(p_b(x\vee y))\\
 &=p_a(p_b(x)\vee p_b(y))\\
 &=p_a(p_b(x))\vee p_a(p_b(y))\\
 &=p_{ab}(x)\vee p_{ab}(y).
\end{align*}
so $ab\in D(K)$. If $a,b\in R(K)$, then $p_a$ and $p_b$ are left adjoints; left adjoints are closed under composition, so $p_{ab}=p_a\circ p_b$ is again a left adjoint, and $ab\in R(K)$.
\end{proof}

For $a\in K$, let $i_a:\downarrow a\hookrightarrow K$ be the inclusion, and define $q_a:K\longrightarrow\downarrow a,\ q_a(x)=ax$.

\begin{proposition}\label{prop:local-meet}
Suppose that $K$ is an ordered incline-type algebra and that $a\in K$. Then the following are equivalent:
\begin{enumerate}
\item $a\in M(K)$;
\item for every $x\in K$, the product $ax$ is the meet $a\wedge x$;
\item $i_a\dashv q_a$.
\end{enumerate}
\end{proposition}

\begin{proof}
(1)$\Longrightarrow$(2) by Lemma~\ref{lem:basic-product}, $ax\le a$ and $ax\le x$, so $ax$ is a lower bound of $a$ and $x$. To know that it is the greatest one, let $z\le a$ and $z\le x$. Then $z\in\downarrow a=\Fix(p_a)$, so $az=z$. Since $z\le x$, by isotonicity of multiplication $az\le ax$, hence $z\le ax$. Therefore $ax=a\wedge x$.

(2)$\Longrightarrow$(3) since $ax\le a$ by Lemma~\ref{lem:basic-product}, the mapping $q_a(x)=ax$ indeed takes values in $\downarrow a$. For $y\in\downarrow a$ and $x\in K$,
\[
 i_a(y)\le x
 \quad\Longleftrightarrow\quad
 y\le x
 \quad\Longleftrightarrow\quad
 y\le a\wedge x=ax=q_a(x),
\]
which is exactly the adjunction $i_a\dashv q_a$.

(3)$\Longrightarrow$(1) we prove $\Fix(p_a)=\downarrow a$. If $x\le a$, then $i_a(x)=x\le x$, and by the adjunction we have $x\le q_a(x)=ax$. On the other hand, by Lemma~\ref{lem:basic-product}, $ax\le x$. Hence $ax=x$, so $x\in\Fix(p_a)$. Conversely, if $ax=x$, then by Lemma~\ref{lem:basic-product}, $x=ax\le a$, so $x\in\downarrow a$. Thus $\Fix(p_a)=\downarrow a$, i.e. $a\in M(K)$.
\end{proof}

\begin{theorem}\label{thm:saturation}
Suppose that $K$ is an ordered incline-type algebra. Then:
\begin{enumerate}
\item
 $E(K)=K\Longleftrightarrow M(K)=K
 \Longleftrightarrow xy=x\wedge y\ (\forall x,y\in K)$.
In this case the underlying join-semilattice is a bounded lattice.
\item $D(K)=K$ if and only if $K$ is an incline.
\item If $K$ is complete, then $R(K)=K$ if and only if $K$ is a commutative integral quantale.
\item $M(K)=R(K)=K$ if and only if $K$ is a Heyting algebra with multiplication equal to meet. If $K$ is complete, these conditions are equivalent to $K$ being a frame.
\end{enumerate}
\end{theorem}

\begin{proof}
(1) Suppose that $E(K)=K$. Then every element is multiplicatively idempotent, i.e. $x^2=x$. For $x,y\in K$, by Lemma~\ref{lem:basic-product}, $xy\le x$ and $xy\le y$, so $xy$ is a lower bound of $x$ and $y$. If $z\le x$ and $z\le y$, by isotonicity $z^2\le xy$, and since $z^2=z$, we get $z\le xy$. Hence $xy=x\wedge y$. By Proposition~\ref{prop:local-meet}, every $a$ satisfies $ax=a\wedge x$, so $a\in M(K)$, and $M(K)=K$. Proposition~\ref{prop:spectra-basic}(2) yields $M(K)\subseteq E(K)$, so the converse $M(K)=K\Longrightarrow E(K)=K$ holds. When $xy=x\wedge y$, $x=x1\le 1$ which shows that $1$ is the greatest element; hence the join-semilattice is a bounded lattice.

(2) $D(K)=K$ means precisely that every multiplier $p_a(x)=ax$ preserves binary joins, i.e. $a(x\vee y)=ax\vee ay\ (\forall a,x,y)$. This is exactly the incline distributivity law, so it is equivalent to $K$ being an incline.

(3) Suppose now that $K$ is complete. If $R(K)=K$, then each $p_a$ has a right adjoint, hence is a left adjoint and preserves arbitrary joins: $a(\bigvee_i x_i)=\bigvee_i ax_i$. This is the quantale distributivity law, so $K$ is a commutative integral quantale. Conversely, if $K$ is a commutative integral quantale, each $p_a$ preserves arbitrary joins; on a complete lattice this forces a right adjoint $r_a(y)=\bigvee\{x:ax\le y\}$, so $a\in R(K)$.

\begin{revblock}
(4) Suppose that $M(K)=K$. Then by (1) multiplication is meet, i.e. $p_a(x)=a\wedge x$. If also $R(K)=K$, let $r_a$ denote the right adjoint of $a\wedge(-)$ and define
\[
 a\to y:=r_a(y).
\]
Then $a\wedge x\le y$ if and only if $x\le a\to y$, so the bounded lattice underlying $K$ is a Heyting algebra. Notice that no completeness assumption is used here; in the complete case the same right adjoint can equivalently be written as $a\to y=\bigvee\{x:a\wedge x\le y\}$. Conversely, if the underlying lattice is Heyting and multiplication is meet, then $a\wedge(-)\dashv(a\to(-))$, so every $p_a$ has a right adjoint and $R(K)=K$; Proposition~\ref{prop:local-meet} also gives $M(K)=K$. When $K$ is complete, a complete Heyting algebra is precisely a frame.
\end{revblock}
\end{proof}

\begin{corollary}
Suppose that $K$ is an ordered incline-type algebra. Then the following are equivalent:
\begin{enumerate}
\item $K$ is a distributive lattice with multiplication equal to meet;
\item $E(K)=K$ and $D(K)=K$;
\item all principal derivations are idempotent and preserve binary joins.
\end{enumerate}
\end{corollary}

\begin{proof}
(1)$\Longleftrightarrow$(2) by Theorem~\ref{thm:saturation}(1), $E(K)=K$ is equivalent to multiplication being meet and to $K$ being a bounded lattice. By part (2), $D(K)=K$ is equivalent to the incline distributivity law, which for meet multiplication says exactly that $K$ is distributive.

(2)$\Longleftrightarrow$(3) by the definitions of the principal spectra, $a\in E(K)$ means $p_a^2=p_a$ and $a\in D(K)$ means that $p_a$ preserves binary joins; hence (2) holds if and only if every principal derivation is idempotent and preserves binary joins.
\end{proof}

The intersection $M(K)\cap R(K)$ has a useful local interpretation.

\begin{proposition}\label{prop:local-adjoint}
Suppose that $K$ is an ordered incline-type algebra and that $a\in M(K)\cap R(K)$. Let $i_a:\downarrow a\hookrightarrow K,\ q_a(x)=ax$, and let $r_a$ be the right adjoint of $p_a$. Then:
\begin{enumerate}
\item $i_a\dashv q_a\dashv r_a|_{\downarrow a}$, where $r_a|_{\downarrow a}:\downarrow a\longrightarrow K$ denotes the restriction of $r_a$ to the domain $\downarrow a$;
\item $p_a$ is an interior operator and $r_a$ is a closure operator;
\item the fixed-point posets $\Fix(p_a)$ and $\Fix(r_a)$ are order-isomorphic.
\end{enumerate}
\end{proposition}

\begin{proof}
(1) The first adjunction $i_a\dashv q_a$ is Proposition~\ref{prop:local-meet}. For $x\in K$ and $y\le a$, since $q_a(x)=ax=p_a(x)$,
\[
 q_a(x)\le y
 \quad\Longleftrightarrow\quad
 p_a(x)\le y
 \quad\Longleftrightarrow\quad
 x\le r_a(y),
\]
so $q_a\dashv r_a|_{\downarrow a}$.

(2) By Lemma~\ref{lem:basic-product}, $p_a(x)=ax\le x$, so $p_a\le\id$. Since $a\in M(K)\subseteq E(K)$ by Proposition~\ref{prop:spectra-basic}(2), we have $a^2=a$, hence $p_a^2=p_{a^2}=p_a$. Thus $p_a$ is an interior operator. Since $p_a^2=p_a$ and $p_a\dashv r_a$, by adjoint uniqueness we have $r_a^2=r_a$; moreover $p_a\le\id$ implies $\id\le r_a$. Hence $r_a$ is a closure operator.

(3) Since $a\in M(K)$, $\Fix(p_a)=\downarrow a$. Consider $r_a:\Fix(p_a)\longrightarrow\Fix(r_a)$ and $p_a:\Fix(r_a)\longrightarrow\Fix(p_a)$. If $y\in\Fix(p_a)$, then $r_a(y)\in\Fix(r_a)$ because $r_a^2(y)=r_a(y)$; similarly, if $z\in\Fix(r_a)$, then $p_a(z)\in\Fix(p_a)$ because $p_a^2(z)=p_a(z)$. The two maps are mutually inverse: for $y\in\Fix(p_a)$, by the adjunction we have $p_a(r_a(y))\le y$, while $y\le r_a(y)$ and monotonicity of $p_a$ yield $y=p_a(y)\le p_a(r_a(y))$; hence $p_a(r_a(y))=y$. Likewise, for $z\in\Fix(r_a)$, $p_a(z)\le z$ and monotonicity of $r_a$ yield $r_a(p_a(z))\le r_a(z)=z$, while by the adjunction we have $z=r_a(z)\le r_a(p_a(z))$; hence $r_a(p_a(z))=z$. Thus $r_a$ and $p_a$ restrict to mutually inverse order-isomorphisms, so $\downarrow a\cong\Fix(r_a)$.
\end{proof}

\section{The complete derivation lattice and its fibers}

Throughout this section, suppose that $K$ is a commutative integral quantale, hence in particular an ordered incline-type algebra. By Theorem~\ref{thm:quantale-characterization}, arbitrary joins in $\Der(K)$ are computed pointwise. Notice again that the individual derivations in our sense are \emph{not} assumed to preserve arbitrary joins; this is the main distinction from the quantale derivations considered in \cite{XiaoLiu2016}.

\begin{revblock}
A complete-lattice result is already known in the important special case of complete MV-algebras: Zhao, Gan and Yang proved that arbitrary pointwise joins of $(\odot,\vee)$-derivations remain derivations and that the resulting derivation poset is complete \cite{ZhaoGanYang2024}. Theorem~\ref{thm:complete-der} extends this phenomenon to all commutative integral quantales under the weaker Leibniz-only definition and, in addition, gives an explicit interior correction for arbitrary meets. Theorem~\ref{thm:fibers} further yields the adjoint triple and the complete-interval decomposition of the evaluation fibers.
\end{revblock}

\begin{definition}
For $f\in K^K$, define the \emph{derivation interior} 
\[
\partial(f)=\bigvee\{d\in\Der(K):d\le f\},
\]
where the join is pointwise.
\end{definition}

\begin{theorem}\label{thm:complete-der}
Suppose that $K$ is a commutative integral quantale. Then the following hold:
\begin{enumerate}
\item $\partial:K^K\longrightarrow K^K$ is an interior operator with $\Fix(\partial)=\Der(K)$; consequently, $\Der(K)$ is a complete lattice;
\item joins in $\Der(K)$ are computed pointwise;
\item meets in $\Der(K)$ are obtained by applying $\partial$ to pointwise meets.
\end{enumerate}
For every family $\{d_i\}_{i\in I}\subseteq\Der(K)$, these take the concrete form
\begin{equation}\label{eq:der-join}
 \left(\bigvee_{i\in I}^{\Der}d_i\right)(x)=\bigvee_{i\in I}d_i(x),
\end{equation}
and
\begin{equation}\label{eq:der-meet}
 \bigwedge_{i\in I}^{\Der}d_i
 =\partial\left(\bigwedge_{i\in I}^{\pt}d_i\right).
\end{equation}
\end{theorem}

\begin{proof}
(1) By Theorem~\ref{thm:quantale-characterization}, arbitrary pointwise joins of derivations are again derivations. Hence $\partial(f)$, being the pointwise join of all derivations below $f$, is itself a derivation and satisfies $\partial(f)\le f$. If $f\le g$, then every derivation below $f$ is also below $g$, so $\partial(f)\le\partial(g)$. Moreover, $\partial(f)$ is a derivation below $f$, when it belongs to the set whose join defines $\partial(f)$, and therefore $\partial(\partial(f))=\partial(f)$. Thus $\partial$ is an interior operator. Its fixed points are exactly the derivations: if $d\in\Der(K)$, then $d$ is the greatest derivation below itself, so $\partial(d)=d$; conversely, if $\partial(f)=f$, then $f=\partial(f)$ is a derivation.

(2) let $d(x)=\bigvee_i d_i(x)$. By Theorem~\ref{thm:quantale-characterization}, $d\in\Der(K)$. It is clearly an upper bound, and if $h\in\Der(K)$ satisfies $d_i\le h$ for all $i$, then $d(x)=\bigvee_i d_i(x)\le h(x)$, so $d\le h$. Hence $d$ is the least upper bound, proving \eqref{eq:der-join}.

(3) let $f=\bigwedge_i^{\pt}d_i$ be the pointwise meet in $K^K$. This $f$ need not be a derivation. Since $\partial(f)$ is the greatest derivation below $f$, we have $\partial(f)\le f\le d_i$ for every $i$, so $\partial(f)$ is a derivation below all $d_i$. Conversely, if $d\in\Der(K)$ satisfies $d\le d_i$ for all $i$, then $d\le f$ pointwise, and by maximality of $\partial(f)$ below $f$ we have $d\le\partial(f)$. Thus $\partial(f)$ is the greatest derivation below all $d_i$, i.e. the meet in $\Der(K)$, which is \eqref{eq:der-meet}.
\end{proof}

The correction in \eqref{eq:der-meet} is necessary in general. Even for lattice derivations, pointwise meets need not remain derivations; see \cite{GanGuo2022}. Theorem~\ref{thm:complete-der} therefore differs in two ways from the complete-lattice constructions for more restrictive derivation notions: joins of our Leibniz derivations are automatic in the quantale case, whereas meets require an interior correction.

We now describe the fibers of evaluation at the unit.

\begin{definition}
For $a\in K$, let $\Der_a(K)=\{d\in\Der(K):d(1)=a\}$, and define $u_a=\bigvee\Der_a(K)$, where the join is taken in $\Der(K)$, equivalently pointwise.
\end{definition}

\begin{theorem}\label{thm:fibers}
Suppose that $K$ is a commutative integral quantale. Then:
\begin{enumerate}
\item $\eta\dashv\varepsilon\dashv u$, where $\eta(a)=p_a$, $\varepsilon(d)=d(1)$, and $u(a)=u_a$;
\item $\Der_a(K)=[p_a,u_a]_{\Der(K)}$;
\item every fiber $\Der_a(K)$ is a complete interval, hence a complete lattice in its induced order; moreover it is closed under all nonempty joins and meets computed in $\Der(K)$;
\item $\varepsilon:\Der(K)\longrightarrow K$ is a surjective complete-lattice homomorphism;
\item define $d\sim e$ by $d(1)=e(1)$. Then
\[
 \Der(K)/{\sim}\ \cong K
\]
as complete lattices.
\end{enumerate}
\end{theorem}

\begin{proof}
(1) The adjunction $\eta\dashv\varepsilon$ is Theorem~\ref{thm:principal}: $p_a\le d$ if and only if $a\le d(1)$. We prove $\varepsilon\dashv u$. First, $u_a\in\Der_a(K)$: since joins in $\Der(K)$ are pointwise,
\[
 u_a(1)=\Bigl(\bigvee_{d\in\Der_a(K)}d\Bigr)(1)=\bigvee_{d\in\Der_a(K)}d(1)=\bigvee_{d\in\Der_a(K)}a=a.
\]

Suppose that $d(1)\le a$. By Theorem~\ref{thm:quantale-characterization}, $d\vee_{\pt}p_a$ is a derivation, and
\[
 (d\vee_{\pt}p_a)(1)=d(1)\vee p_a(1)=d(1)\vee a=a.
\]
Thus $d\vee_{\pt}p_a\in\Der_a(K)$, so $d\vee_{\pt}p_a\le u_a$, and hence $d\le u_a$. Conversely, if $d\le u_a$, then
\[
 d(1)\le u_a(1)=a.
\]
Therefore $\varepsilon(d)\le a$ if and only if $d\le u_a$, i.e. $\varepsilon\dashv u$.

(2) By Theorem~\ref{thm:principal}, $p_a\le d$ for every $d\in\Der_a(K)$, while $u_a$ is the largest element of the fiber. Conversely, if $p_a\le d\le u_a$, when it takes value 1 then we have $a=p_a(1)\le d(1)\le u_a(1)=a$, so $d(1)=a$, i.e. $d\in\Der_a(K)$.

(3) Every interval in a complete lattice is complete in its induced order, with bottom $p_a$ and top $u_a$. Since $\varepsilon$ is simultaneously a right and a left adjoint, it preserves arbitrary joins and meets. Hence for any nonempty family $\{d_i\}_{i\in I}\subseteq\Der_a(K)$,
\[
 \varepsilon\Bigl(\bigvee_i d_i\Bigr)=\bigvee_i\varepsilon(d_i)=\bigvee_i a=a,
\quad
 \varepsilon\Bigl(\bigwedge_i d_i\Bigr)=\bigwedge_i\varepsilon(d_i)=\bigwedge_i a=a,
\]
so the nonempty join or meet remains in the fiber.

(4) The same adjoint argument shows that $\varepsilon$ is a complete-lattice homomorphism; $\varepsilon(p_a)=p_a(1)=a$ shows surjectivity.

(5) The relation $\sim$ is the kernel congruence of the surjective complete-lattice homomorphism $\varepsilon$, so the quotient is isomorphic to $K$.
\end{proof}

\begin{corollary}\label{cor:meet-core}
For every family $\{d_i\}_{i\in I}\subseteq\Der(K)$,
\[
 p_{\bigwedge_{i\in I} d_i(1)}
 \le
 \bigwedge_{i\in I}^{\Der}d_i
 \le
 \bigwedge_{i\in I}^{\pt}d_i.
\]
Thus the left-hand term is the principal core of the derivation meet.
\end{corollary}

\begin{proof}
Let $m=\bigwedge_i^{\Der}d_i$. Since $\varepsilon$ preserves meets, $m(1)=\varepsilon(m)=\bigwedge_i\varepsilon(d_i)=\bigwedge_i d_i(1)$.
By Theorem~\ref{thm:principal}(4), $p_{m(1)}\le m$, hence the left inequality. For the right inequality, $m$ is a pointwise lower bound of the $d_i$ in $K^K$, hence $m\le\bigwedge_i^{\pt}d_i$.
\end{proof}

\section{\rev{t-Norms: reconstruction, Archimedean classification, and information loss}}

We now specialize to $K_T=([0,1],\max,T,0,1)$. For a t-norm $T$, write $\Der(T)$, $\Der_a(T)$, and $\PDer(T)$ instead of $\Der(K_T)$, $\Der_a(K_T)$, and $\PDer(K_T)$, and put $p_a^T(x)=T(a,x)$.

\begin{revblock}
\begin{definition}\label{def:fibered-poset}
A \emph{fibered poset over $[0,1]$} is a pair $(P,\varepsilon_P)$ consisting of a poset $P$ and an isotone map $\varepsilon_P:P\to[0,1]$. An isomorphism over $[0,1]$ from $(P,\varepsilon_P)$ to $(Q,\varepsilon_Q)$ is an order isomorphism $\Phi:P\to Q$ satisfying $\varepsilon_Q\circ\Phi=\varepsilon_P$. More generally, if $h:[0,1]\to[0,1]$ is an increasing bijection, we say that $\Phi$ \emph{covers $h$} when $\varepsilon_Q\circ\Phi=h\circ\varepsilon_P$.
\end{definition}

For a t-norm $T$, set $\mathfrak F(T)=(\Der(T),\varepsilon_T)$, where $\varepsilon_T(d)=d(1)$, and set $\mathfrak P(T)=(\PDer(T),\le,\circ,p_1^T)$. The first invariant records order and evaluation fibers; the second additionally remembers composition of the canonical fiber minima.
\end{revblock}

\subsection{Principal reconstruction}

\begin{theorem}\label{thm:tnorm-reconstruction}
Suppose that $T$ is a t-norm. Then:
\begin{enumerate}
\item $p_a^T=\min\Der_a(T)$ for every $a\in[0,1]$;
\item $p_a^T\circ p_b^T=p_{T(a,b)}^T$;
\item $([0,1],\le,T,1)\cong(\PDer(T),\le,\circ,p_1^T)$;
\item the t-norm is reconstructed concretely by $T(a,b)=\bigl(p_a^T\circ p_b^T\bigr)(1)$.
\end{enumerate}
\end{theorem}

\begin{proof}
(1) By Theorem~\ref{thm:principal}(4), $p_a^T\le d$ holds whenever $a\le d(1)$. For $d\in\Der_a(T)$ we have $d(1)=a$, so $p_a^T\le d$; since $p_a^T(1)=a$, it lies in $\Der_a(T)$ and is therefore its least element. This argument uses only the ordered incline-type axioms and hence applies to every t-norm $T$.

(2) By Theorem~\ref{thm:principal}(2) we have $p_a\circ p_b=p_{ab}$. In $K_T$ the monoid operation is $T$, so explicitly $(p_a^T\circ p_b^T)(x)=T(a,T(b,x))=T(T(a,b),x)=p_{T(a,b)}^T(x)$, where the middle equality uses associativity of $T$.

(3) The map $\Phi(a)=p_a^T$ is order-preserving since $a\le b$ implies $T(a,x)\le T(b,x)$. It preserves the monoid operation by (2), sends $1$ to the unit $p_1^T$, and is injective because $\Phi(a)=\Phi(b)$ evaluated at $1$ yields $a=T(a,1)=T(b,1)=b$. It is surjective by definition of principal derivations. Hence $\Phi$ is an order-monoid isomorphism.

(4) Evaluate (2) at $x=1$: $(p_a^T\circ p_b^T)(1)=p_{T(a,b)}^T(1)=T(T(a,b),1)=T(a,b)$, using that $1$ is the unit of $T$.
\end{proof}

\begin{revblock}
\begin{proposition}[Universal endpoint-modified principal derivations]\label{prop:endpoint-universal}
Let $T$ be any t-norm and let $0\le a\le c\le1$. Define
\[
 d_{a,c}^T(x)=
 \begin{cases}
 T(c,x),&0\le x<1,\\
 a,&x=1.
 \end{cases}
\]
Then $d_{a,c}^T\in\Der(T)$.
\end{proposition}

\begin{proof}
If $x,y<1$, then $T(x,y)\le\min\{x,y\}<1$, and associativity and commutativity give
\[
 d_{a,c}^T(T(x,y))=T(c,T(x,y))=T(T(c,x),y)=T(x,T(c,y)).
\]
Thus the two Leibniz terms are equal to the left-hand side. If $x=1$ and $y<1$, then $T(a,y)\le T(c,y)=d_{a,c}^T(y)$, so the maximum of the two Leibniz terms is $d_{a,c}^T(y)$. The case $y=1$ is symmetric, and for $x=y=1$ both sides equal $a$.
\end{proof}
\end{revblock}

\subsection{Product derivations}

Let
\[
 \TP(x,y)=xy.
\]

\begin{proposition}\label{prop:product-classification}
A map $d:[0,1]\longrightarrow[0,1]$ is a derivation of $\TP$ if and only if there exist $0\le a\le c\le1$ such that
\[
 d(x)=
 \begin{cases}
 cx,&0\le x<1,\\
 a,&x=1.
 \end{cases}
\]
We denote this derivation by $d_{a,c}^{\mathrm P}$.
\end{proposition}

\begin{proof}
($\Longrightarrow$) let $d$ be a derivation. For $x\in(0,1)$ define $f(x)=\frac{d(x)}{x}$, so that $d(x)=xf(x)$. For $x,y\in(0,1)$ the derivation identity yields $d(xy)=\max\{d(x)y,xd(y)\}$.
Substituting $d(x)=xf(x)$ and $d(y)=yf(y)$ on both sides, $xy\,f(xy)=\max\{xyf(x),xyf(y)\}=xy\max\{f(x),f(y)\}$.
Since $xy>0$, canceling $xy$ yields $f(xy)=\max\{f(x),f(y)\}$.
Set $g(t)=f(e^{-t})$ for $t>0$. Then $g(s+t)=f(e^{-(s+t)})=f(e^{-s}e^{-t})=\max\{f(e^{-s}),f(e^{-t})\}=\max\{g(s),g(t)\}$.
Given $s,t>0$, choose $n$ with $nt>s$. By repeated application of the identity, $g(2t)=\max\{g(t),g(t)\}=g(t)$, and inductively $g(nt)=g(t)$. Hence $g(t)=g(nt)=g(s+(nt-s))=\max\{g(s),g(nt-s)\}\ge g(s)$.
Interchanging $s$ and $t$ yields equality. Hence $g$, and therefore $f$, is constant; say $f(x)=c$. Thus $d(x)=cx$ for $0<x<1$, while by Lemma~\ref{lem:d-basic}, $d(0)=0$. Since $cx=d(x)\le1$ for every $0<x<1$, necessarily $c\le1$. Let $a=d(1)$. By Lemma~\ref{lem:d-basic}, $ax\le d(x)=cx\ (0<x<1)$, so $a\le c$.

($\Longleftarrow$) fix $0\le a\le c\le1$ and define $d$ as above. If $x,y<1$, then $d(xy)=cxy=\max\{cxy,cxy\}$.
If $x=1$ and $y<1$, then $\max\{d(1)y,1\cdot d(y)\}=\max\{ay,cy\}=cy=d(y)$, and the symmetric case is identical. If $x=y=1$, both sides equal $a$. Thus $d$ is a derivation.
\end{proof}

\subsection{{\L}ukasiewicz derivations}

Let
\[
 \TL(x,y)=\max\{0,x+y-1\}.
\]
The standard MV-algebra on $[0,1]$ has $x\odot y=\TL(x,y)$, so the derivations considered in this subsection are precisely $(\odot,\vee)$-derivations in the sense of Zhao, Gan and Yang \cite{ZhaoGanYang2024}. Their work treats general MV-algebras and, in particular, enumerates the derivations on finite MV-chains. The following proposition provides an explicit classification for the standard continuum chain; this classification is then compared with the corresponding one for the Product t-norm.

\begin{proposition}\label{prop:luk-classification}
A map $d:[0,1]\longrightarrow[0,1]$ is a derivation of $\TL$ if and only if there exist $0\le a\le c\le1$ such that
\[
 d(x)=
 \begin{cases}
 \TL(c,x),&0\le x<1,\\
 a,&x=1.
 \end{cases}
\]
We denote this derivation by $d_{a,c}^{\mathrm L}$.
\end{proposition}

\begin{proof}
($\Longrightarrow$) let $d$ be a derivation and put
\[
 F(\alpha)=1-d(1-\alpha),\quad 0\le \alpha\le1.
\]
If $\alpha+\beta\le1$, then
\[
 \TL(1-\alpha,1-\beta)=\max\{0,(1-\alpha)+(1-\beta)-1\}=1-\alpha-\beta.
\]
Substituting $x=1-\alpha$, $y=1-\beta$ into the derivation identity and rewriting in terms of $F$ yields
\begin{equation}\label{eq:F}
 F(\alpha+\beta)=\min\{1,F(\alpha)+\beta,\alpha+F(\beta)\}.
\end{equation}
Set $G(\alpha)=F(\alpha)-\alpha$. Subtracting $\alpha+\beta$ from \eqref{eq:F} yields
\begin{equation}\label{eq:G}
 G(\alpha+\beta)=\min\{G(\alpha),G(\beta),1-\alpha-\beta\}
 \ (\alpha+\beta\le1).
\end{equation}

Taking $\beta=1-\alpha$ in \eqref{eq:G} yields $G(1)=\min\{G(\alpha),G(1-\alpha),0\}$.
By Lemma~\ref{lem:d-basic}, $d(0)=0$, so $F(1)=1-d(0)=1$ and $G(1)=F(1)-1=0$. Hence $0=\min\{G(\alpha),G(1-\alpha),0\}$, which forces $G(\alpha)\ge0$. If $0<\alpha<\beta<1$, write $\beta=\alpha+(\beta-\alpha)$ in \eqref{eq:G}; then $G(\beta)\le G(\alpha)$, so $G$ is nonincreasing on $(0,1)$.

Fix $\alpha\in(0,1)$ and write $\alpha=n(\alpha/n)$. Iterating \eqref{eq:G} gives
\[
 G(\alpha)=\min\{G(\alpha/n),1-\alpha\}.
\]
Let
\[
 q=\sup_{0<t\le1}G(t)=\lim_{t\downarrow0}G(t).
\]
The last equality uses only the monotonicity of $G$ on $(0,1)$: as $t\downarrow0$, the values $G(t)$ increase to their supremum. Hence $G(\alpha/n)\to q$, and passing to the limit gives
\[
 G(\alpha)=\min\{q,1-\alpha\}.
\]

Setting $\alpha=0$ in \eqref{eq:G} yields $G(\beta)=\min\{G(0),G(\beta),1-\beta\}\le G(0)$ for all $\beta\in[0,1]$, hence $q=\sup_{0<t\le1}G(t)\le G(0)$. Since $G(0)=F(0)=1-d(1)$, if $a=d(1)$ then $q\le1-a$. Set $c=1-q$; thus $a\le c\le1$.

For $x<1$, put $\alpha=1-x>0$. Then
\[
 G(1-x)=\min\{q,x\}=\min\{1-c,x\},
\]
and
\[
 F(1-x)=(1-x)+G(1-x)=1-x+\min\{1-c,x\}.
\]
Therefore
\[
 d(x)=1-F(1-x)=x-\min\{1-c,x\}=\max\{0,c+x-1\}=\TL(c,x).
\]

($\Longleftarrow$) fix $0\le a\le c\le1$ and define $d$ by the displayed formula. If $x,y<1$, then $\TL(x,y)<1$, and associativity and commutativity of $\TL$ yield
\begin{align*}
 d(\TL(x,y))
 &=\TL(c,\TL(x,y))\\
 &=\TL(\TL(c,x),y)
 =\TL(x,\TL(c,y)).
\end{align*}
Hence the maximum of the two Leibniz terms equals $d(\TL(x,y))$. If $x=1$ and $y<1$, then $\max\{\TL(a,y),\TL(1,d(y))\}=\max\{\TL(a,y),\TL(c,y)\}=\TL(c,y)=d(y)$, because $a\le c$. The case $y=1$ is symmetric, and for $x=y=1$ both sides equal $a$. Thus $d$ is a derivation.
\end{proof}

\begin{revblock}
\subsection{Continuous Archimedean t-norms}

A continuous Archimedean t-norm is called \emph{strict} when it has no nonzero nilpotent elements, and \emph{nilpotent} when every $x\in(0,1)$ has a finite power equal to $0$. The classical representation theorem states that every continuous Archimedean t-norm is isomorphic, as an ordered monoid, either to $\TP$ in the strict case or to $\TL$ in the nilpotent case \cite{KlementMesiarPap2000,Jenei2002,MostertShields1957}. Equivalently, there is an increasing bijection $h$ satisfying $h(T(x,y))=S(h(x),h(y))$ for the corresponding standard model $S$. The following transport lemma turns the two model classifications into a classification of the entire Archimedean class.

\begin{lemma}[Transport of derivations]\label{lem:transport}
Let $T$ and $S$ be t-norms, and let $h:[0,1]\to[0,1]$ be an increasing bijection satisfying $h(T(x,y))=S(h(x),h(y))$. Then
\[
 \Phi_h:\Der(T)\longrightarrow\Der(S),\qquad \Phi_h(d)=h\circ d\circ h^{-1},
\]
is an order isomorphism covering $h$. Moreover, $\Phi_h(p_c^T)=p_{h(c)}^S$.
\end{lemma}

\begin{proof}
An increasing bijection of $[0,1]$ preserves $0$, $1$, and binary maxima. If $u=h(x)$ and $v=h(y)$, then
\begin{align*}
 \Phi_h(d)(S(u,v))
 &=h(d(T(x,y)))\\
 &=h\bigl(T(d(x),y)\vee T(x,d(y))\bigr)\\
 &=S(\Phi_h(d)(u),v)\vee S(u,\Phi_h(d)(v)).
\end{align*}
Thus $\Phi_h(d)$ is a derivation. Its inverse is $\Phi_{h^{-1}}$, and pointwise order is preserved in both directions. Finally, $\varepsilon_S(\Phi_h(d))=h(d(1))=h(\varepsilon_T(d))$, and the principal identity follows directly from the t-norm isomorphism.
\end{proof}

\begin{theorem}\label{thm:arch-classification}
Let $T$ be a continuous Archimedean t-norm and let $d:[0,1]\to[0,1]$. Then $d$ is a derivation of $T$ if and only if there exist unique $a,c\in[0,1]$ with $a\le c$ such that
\[
 d=d_{a,c}^T.
\]
\end{theorem}

\begin{proof}
Choose an increasing t-norm isomorphism $h$ from $T$ onto either $\TP$ or $\TL$. By Lemma~\ref{lem:transport}, $\widetilde d=h\circ d\circ h^{-1}$ is a derivation of the corresponding standard model. Propositions~\ref{prop:product-classification} and \ref{prop:luk-classification} give parameters $A\le C$ such that $\widetilde d(u)=S(C,u)$ for $u<1$ and $\widetilde d(1)=A$. Put $a=h^{-1}(A)$ and $c=h^{-1}(C)$. Then $d(x)=T(c,x)$ for $x<1$ and $d(1)=a$. The converse is Proposition~\ref{prop:endpoint-universal}. For uniqueness, $a=d(1)$ is forced. If both $c$ and $c'$ represent the same $d$, then $T(c,x)=T(c',x)$ for every $x<1$; letting $x\uparrow1$ and using continuity of $T$ gives $c=c'$.
\end{proof}

\begin{theorem}[Universal Archimedean fibered-poset no-go theorem]\label{thm:universal-no-go}
Let $\Delta=\{(a,c)\in[0,1]^2:a\le c\}$ with the coordinatewise order, and let $\pi_1(a,c)=a$. For every continuous Archimedean t-norm $T$, the map $d_{a,c}^T\mapsto(a,c)$ is an isomorphism over $[0,1]$ from $\mathfrak F(T)$ onto $(\Delta,\pi_1)$. Consequently, any two continuous Archimedean t-norms have isomorphic fibered derivation posets.

Nevertheless, the fibered derivation poset does not determine whether $T$ is strict or nilpotent. By contrast, the ordered composition monoid $\mathfrak P(T)$ determines the complete t-norm through Theorem~\ref{thm:tnorm-reconstruction}.
\end{theorem}

\begin{proof}
Theorem~\ref{thm:arch-classification} gives a bijection $\Der(T)\to\Delta$, and evaluation at $1$ is the first coordinate. Moreover,
\[
 d_{a,c}^T\le d_{b,e}^T\quad\Longleftrightarrow\quad a\le b\text{ and }c\le e.
\]
The reverse implication follows from isotonicity of $T$. For the forward implication, evaluation at $1$ gives $a\le b$, while $T(c,x)\le T(e,x)$ for every $x<1$; letting $x\uparrow1$ and using continuity yields $c\le e$. Hence the displayed map is an isomorphism over $[0,1]$.

To see the information loss, compare one strict and one nilpotent continuous Archimedean t-norm. Their fibered derivation posets are isomorphic, whereas their ordered monoids are not: the nilpotent type has nonzero nilpotent elements and the strict type does not.
\end{proof}

\begin{remark}[Where the missing information resides]\label{rem:missing-composition}
Under the universal identification with $\Delta$, the minimum of the fiber over $a$ is always
\[
 p_a^T=d_{a,a}^T\longleftrightarrow(a,a).
\]
Hence the fibered poset already determines the distinguished family of principal minima, and this enriched ordered data is still independent of the continuous Archimedean t-norm. The t-norm-dependent information first appears when composition is retained, since
\[
 p_a^T\circ p_b^T=p_{T(a,b)}^T
 \qquad\text{corresponds to}\qquad
 (a,a)\circ_T(b,b)=(T(a,b),T(a,b)).
\]
Thus Theorem~\ref{thm:universal-no-go} isolates composition, rather than the order or the fiber minima themselves, as the missing reconstruction datum.
\end{remark}

\begin{corollary}[Explicit Archimedean fibers and lattice operations]\label{cor:arch-fibers}
Let $T$ be a continuous Archimedean t-norm. Then, for every $a\in[0,1]$,
\[
 \Der_a(T)=\{d_{a,c}^T:c\in[a,1]\},
\]
and the map $c\mapsto d_{a,c}^T$ is an order isomorphism $[a,1]\cong\Der_a(T)$. In particular,
\[
 p_a^T=d_{a,a}^T,
 \qquad
 u_a=d_{a,1}^T.
\]
Moreover, for every family $\{(a_i,c_i)\}_{i\in I}\subseteq\Delta$,
\[
 \bigvee_{i\in I}^{\Der}d_{a_i,c_i}^T
 =d_{\bigvee_i a_i,\,\bigvee_i c_i}^T,
 \qquad
 \bigwedge_{i\in I}^{\Der}d_{a_i,c_i}^T
 =d_{\bigwedge_i a_i,\,\bigwedge_i c_i}^T.
\]
Thus the order isomorphism in Theorem~\ref{thm:universal-no-go} is a complete-lattice isomorphism after the evaluation map is forgotten.
\end{corollary}

\begin{proof}
Theorem~\ref{thm:arch-classification} gives the description of the fiber, and the order calculation in the proof of Theorem~\ref{thm:universal-no-go} shows that $c\mapsto d_{a,c}^T$ is an order isomorphism. Since $p_a^T(x)=T(a,x)$ also at $x=1$, one has $p_a^T=d_{a,a}^T$; the largest permissible second parameter is $1$, so the fiber maximum is $u_a=d_{a,1}^T$.

The subset $\Delta\subseteq[0,1]^2$ is closed under coordinatewise arbitrary joins and meets. Hence the order isomorphism $\Der(T)\cong\Delta$ transports those operations and yields the displayed formulas. For the empty family these formulas give the ambient bottom $d_{0,0}^T$ and ambient top $d_{1,1}^T$. By contrast, the internal bottom and top of the fixed fiber $\Der_a(T)$ are $d_{a,a}^T$ and $d_{a,1}^T$, respectively, in agreement with Theorem~\ref{thm:fibers}.
\end{proof}
\end{revblock}

\section{Cut derivations and the idempotent skeleton}\label{sec:cut-spectrum}

Suppose now that $T$ is continuous. Let $E(T)=\{e\in[0,1]:T(e,e)=e\}$ be its idempotent set. For $e\in[0,1]$, define the canonical cut map
\[
 \kappa_e(x)=
 \begin{cases}
 x,&x<e,\\
 0,&x\ge e.
 \end{cases}
\]
This family is generally nonprincipal, although some endpoint cases may be principal.

\begin{revblock}
Threshold maps of the same form already occur in Liu and Liu's study of linearly ordered EBL-algebras: their Proposition~4.1 gives a sufficient condition for the map $x\mapsto x$ below a threshold and $0$ at or above it to be a derivation \cite{LiuLiu2024}. For a t-norm and threshold $e$, that condition is $T(x,y)\ge e$ whenever $x,y\ge e$. By isotonicity this is equivalent to $T(e,e)\ge e$, and since every t-norm satisfies $T(e,e)\le e$, it is equivalent to the idempotence of $e$. Thus the earlier sufficient criterion yields the forward implication below. The additional conclusion established here is the converse under continuity; together the two directions make the cut family an exact detector of the whole idempotent skeleton. This also refines earlier derivation-based characterizations of G\"odel algebras \cite{TorkzadehAbbasian2013,XinFengYang2016}.
\end{revblock}

We use the standard property that if $T$ is continuous and $e$ is idempotent, then
\begin{equation}\label{eq:idempotent-local-unit}
 T(x,e)=x\ (0\le x\le e).
\end{equation}
It follows, for example, from the ordinal-sum representation of continuous t-norms; see \cite{KlementMesiarPap2000,Jenei2002,MostertShields1957}.

\begin{revblock}
\begin{theorem}[Cut-spectrum characterization]\label{thm:cut-spectrum}
Let $T$ be a continuous t-norm. For every $e\in[0,1]$,
\[
 e\in E(T)
 \quad\Longleftrightarrow\quad
 \kappa_e\in\Der_0(T).
\]
Consequently,
\[
 E(T)=\{e\in[0,1]:\kappa_e\in\Der_0(T)\}.
\]
\end{theorem}

\begin{proof}
Suppose first that $e$ is idempotent. If $x,y<e$, then $T(x,y)\le\min\{x,y\}<e$, and hence
\[
 \kappa_e(T(x,y))=T(x,y)
 =T(\kappa_e(x),y)\vee T(x,\kappa_e(y)).
\]
If $x<e\le y$, then isotonicity and \eqref{eq:idempotent-local-unit} give
\[
 x=T(x,e)\le T(x,y)\le x,
\]
so $T(x,y)=x$ and the derivation identity follows; the case $y<e\le x$ is symmetric. If $e\le x,y$, then $T(x,y)\ge T(e,e)=e$, so both sides of the derivation identity are $0$. Thus $\kappa_e\in\Der(T)$, and $\kappa_e(1)=0$.

Conversely, suppose that $e$ is not idempotent. Then $T(e,e)<e$. The diagonal map $\delta(x)=T(x,x)$ is continuous, with $\delta(e)<e<\delta(1)=1$. Choose $r$ such that $\delta(e)<r<e$. The intermediate value theorem gives an $x>e$ with $T(x,x)=r$. Then $\kappa_e(x)=0$ but
\[
 \kappa_e(T(x,x))=\kappa_e(r)=r>0,
\]
whereas the derivation identity would make this value equal to
\[
 T(0,x)\vee T(x,0)=0,
\]
a contradiction.
\end{proof}
\end{revblock}

\begin{corollary}
A continuous t-norm $T$ is Archimedean if and only if
\[
 \kappa_e\notin\Der_0(T)
 \quad\text{for every }e\in(0,1).
\]
\end{corollary}

\begin{proof}
A continuous t-norm is Archimedean if and only if its only idempotents are $0$ and $1$. Apply Theorem~\ref{thm:cut-spectrum}.
\end{proof}

\begin{corollary}
Suppose that $T$ is a continuous t-norm. Then the following are equivalent:
\begin{enumerate}
\item $T=\TG$, where $\TG(x,y)=\min\{x,y\}$;
\item $\kappa_e\in\Der_0(T)$ for every $e\in[0,1]$.
\end{enumerate}
\end{corollary}

\begin{proof}
\begin{revblock}
If $T=\TG$, every point is idempotent, so Theorem~\ref{thm:cut-spectrum} gives $\kappa_e\in\Der_0(T)$ for every $e$. Conversely, if all the cut maps belong to $\Der_0(T)$, then Theorem~\ref{thm:cut-spectrum} gives $T(e,e)=e$ for every $e$. For $x\le y$,
\[
 x=T(x,x)\le T(x,y)\le x,
\]
so $T(x,y)=x=\min\{x,y\}$.
\end{revblock}
\end{proof}

\begin{revblock}
\section{Direct operator detection of the ordinal-sum decomposition}

Let $T$ be a continuous t-norm. The set $E(T)$ is closed, and every connected component $(\alpha,\beta)$ of $[0,1]\setminus E(T)$ determines an Archimedean ordinal-sum component on $[\alpha,\beta]$, of either Product or {\L}ukasiewicz type \cite{MostertShields1957,KlementMesiarPap2000,Jenei2002}. The finite-power distinction between strict and nilpotent continuous t-norms is classical; see Walker and Walker \cite{WalkerWalker2002}. Our purpose is to translate this distinction into a finite-composition criterion on principal derivations.

For $c\in[0,1]$, define $c^{(1)}=c$ and $c^{(n+1)}=T(c^{(n)},c)$.

\begin{lemma}\label{lem:principal-powers}
For every $c\in[0,1]$ and $n\ge1$, $(p_c^T)^{\circ n}=p_{c^{(n)}}^T$. Hence $(p_c^T)^{\circ n}=p_\alpha^T$ if and only if $c^{(n)}=\alpha$.
\end{lemma}

\begin{proof}
The first identity follows inductively from $p_a^T\circ p_b^T=p_{T(a,b)}^T$; the second follows from injectivity of the principal representation.
\end{proof}

\begin{theorem}[Principal-derivation form of the strict/nilpotent power criterion]\label{thm:block-type}
Let $(\alpha,\beta)$ be a connected component of $[0,1]\setminus E(T)$. Then:
\begin{enumerate}
\item the component is {\L}ukasiewicz-type if and only if there exist $c\in(\alpha,\beta)$ and $n\ge2$ such that $(p_c^T)^{\circ n}=p_\alpha^T$; equivalently, every $c\in(\alpha,\beta)$ has such a finite power;
\item the component is Product-type if and only if $(p_c^T)^{\circ n}>p_\alpha^T$ for every $c\in(\alpha,\beta)$ and every $n\ge1$; equivalently, it is enough that one interior $c$ never reaches $p_\alpha^T$ by a finite composition power.
\end{enumerate}
\end{theorem}

\begin{proof}
Let $S$ be the standard t-norm associated with the component, and let $h:[\alpha,\beta]\to[0,1]$ be the increasing isomorphism satisfying $h(\alpha)=0$, $h(\beta)=1$, and
\[
 h(T(x,y))=S(h(x),h(y))
 \qquad(x,y\in[\alpha,\beta]).
\]
Induction gives $h(c^{(n)})=h(c)^{(n)}_S$ for every $n$.\par
If the component is of {\L}ukasiewicz type, every $u\in(0,1)$ has a finite {\L}ukasiewicz power equal to $0$. Hence every $c\in(\alpha,\beta)$ has a finite $T$-power equal to $\alpha$.\par
If the component is of Product type, every finite Product power of $u\in(0,1)$ is strictly positive. Hence $c^{(n)}>\alpha$ for every finite $n$. Lemma~\ref{lem:principal-powers} translates these alternatives into the stated composition properties. Since every continuous Archimedean component is of exactly one of the two types, the existential and universal formulations are equivalent.
\end{proof}

\begin{theorem}[Operator recovery of the Mostert--Shields invariants]\label{thm:MS-detection}
Let $T$ be a continuous t-norm. Within the concrete representation $\Der(T)\subseteq[0,1]^{[0,1]}$, the following data are recovered in succession:
\begin{enumerate}
\item the minimum of the evaluation fiber $\Der_a(T)$ is the principal derivation $p_a^T$;
\item the cut spectrum recovers the idempotent set through
\[
 E(T)=\{e\in[0,1]:\kappa_e\in\Der_0(T)\};
\]
\item once $E(T)$ is known, the Archimedean blocks are the connected components of $[0,1]\setminus E(T)$;
\item the type of each block is recovered from finite composition powers of its interior principal derivations, as in Theorem~\ref{thm:block-type}.
\end{enumerate}
Thus principal composition reconstructs the complete t-norm, whereas the cut spectrum and finite principal powers expose its ordinal-sum invariants directly.
\end{theorem}

\begin{proof}
Item (1) is Theorem~\ref{thm:tnorm-reconstruction}(1), and item (2) is Theorem~\ref{thm:cut-spectrum}. Item (3) is the classical ordinal-sum theorem for continuous t-norms applied to the recovered closed set $E(T)$, while item (4) is Theorem~\ref{thm:block-type}. The final comparison follows from Theorem~\ref{thm:tnorm-reconstruction}(2)--(4).
\end{proof}

\begin{example}[A mixed Product--{\L}ukasiewicz ordinal sum]\label{ex:mixed-ordinal-sum}
Define
\[
T(x,y)=
\begin{cases}
2xy,&x,y\in[0,\frac12],\\
\max\{\frac12,x+y-1\},&x,y\in[\frac12,1],\\
\min\{x,y\},&\text{otherwise}.
\end{cases}
\]
On the lower scaled-Product block the only idempotents are $0$ and $\frac12$, while on the upper scaled-{\L}ukasiewicz block they are $\frac12$ and $1$; hence $E(T)=\{0,\frac12,1\}$. Therefore Theorem~\ref{thm:cut-spectrum} identifies the two nontrivial blocks exactly from the cut derivations. For $c=\frac14$, one has $c^{(n)}=\frac12(\frac12)^n>0$ for all $n$, so the lower block is Product-type. For $c=\frac34$, $T(c,c)=\frac12$, so $(p_{3/4}^T)^{\circ2}=p_{1/2}^T$, detecting the {\L}ukasiewicz type of the upper block.
\end{example}
\end{revblock}

\section{Concluding remarks}

\begin{revblock}
The central point of the paper is that derivations become more informative when distributivity and residuation are not built into the ambient algebra from the outset. Principal derivations exist at the weak ordered incline-type level and form a canonical Galois-inserted copy of the multiplicative ordered algebra. Their pointwise-join behavior detects incline and quantale distributivity, their adjoints detect residuation, and their idempotence together with fixed-point behavior detects meet structure. In the complete quantale setting, the full space of Leibniz derivations is a complete lattice with an interior correction for meets and a complete interval decomposition over the underlying quantale.

For t-norms, endpoint-modified principal maps provide a universal construction, and Theorem~\ref{thm:arch-classification} proves that they exhaust all derivations in the continuous Archimedean case. Consequently, the complete fibered derivation poset is universal throughout that class and forgets even the strict/nilpotent type. The ordered composition monoid of principal fiber minima nevertheless reconstructs the multiplication. For general continuous t-norms, cut derivations recover the idempotent skeleton and finite principal composition powers identify the Archimedean component types.

Several questions remain open. First, one may seek a complete description of $\Der_a(T)$ from the ordinal-sum data of a non-Archimedean continuous t-norm and determine how derivations interact across distinct components. Second, noncommutative ordered incline-type structures should produce separate left and right principal spectra. Third, it is natural to ask when $\Der(K)$ is itself a frame or admits a multiplication compatible with its complete-lattice structure. Finally, analogous operator probes for left-continuous and other noncontinuous t-norms remain to be developed.
\end{revblock}


\end{document}